\documentclass[12pt]{amsart}

\usepackage{a4wide}

\usepackage{amssymb,amsmath}

\usepackage{hyperref}

\newtheorem{theorem}{Theorem}
\newtheorem{lemma}[theorem]{Lemma}
\newtheorem{corollary}[theorem]{Corollary}
\newtheorem{proposition}[theorem]{Proposition}

\newtheorem{examples}[theorem]{Examples}
\newtheorem{conjecture}[theorem]{Conjecture}

\newcommand\LF{{\rm F}}

\newcommand{\R}{\ensuremath{\mathbb{R}}}
\newcommand{\C}{\ensuremath{\mathbb{C}}}
\renewcommand{\H}{\ensuremath{\mathbb{H}}}
\newcommand{\Ca}{\ensuremath{\mathbb{O}}}
\renewcommand\SS{\ensuremath{{\sf S}}}

\newcommand{\g}[1]{\ensuremath{\mathfrak{#1}}}

\newcommand{\Isom}{\ensuremath{\mathop{\rm Isom}\nolimits}}

\newcommand{\Ad}{\ensuremath{\mathop{\rm Ad}\nolimits}}

\newcommand{\rank}{\ensuremath{\mathop{\rm rank}\nolimits}}

\newcommand\SO{{\rm SO}}
\newcommand\SUxU[2]{{\rm S(U(}#1)\times{\rm U(}#2{\rm))}}
\newcommand\SU{{\rm SU}}

\newcommand\Spin{{\rm Spin}}
\newcommand\Sp{{\rm Sp}}
\newcommand\U{{\rm U}}

\renewcommand\P{{\sf P}}
\newcommand\HH{{\sf H}}

\begin{document}
\title[Polar actions]{Polar actions with a fixed point}

\begin{abstract}
We prove a criterion for an isometric action of a Lie group on a Riemannian
manifold to be polar. From this criterion, it follows that an action with a
fixed point is polar if and only if the slice representation at the fixed point
is polar and the section is the tangent space of an embedded totally geodesic
submanifold. We apply this to obtain a classification of polar actions with a
fixed point on symmetric spaces.
\end{abstract}

\author[J. C. D\'\i{}az-Ramos]{Jos\'{e} Carlos D\'\i{}az-Ramos}
\address{Department of Geometry and Topology,
University of Santiago de Compostela, Spain.} \email{josecarlos.diaz@usc.es}

\author[A. Kollross]{Andreas Kollross}
\address{Institut f\"{u}r Mathematik,
Universit\"{a}t Augsburg, Universit\"{a}tsstra{\ss}e 14, 86135 Augsburg, Germany}
\email{kollross@math.uni-augsburg.de}

\thanks{The first author has been supported by a Marie-Curie
Fellowship (PERG04-GA-2008-239162) and by projects MTM2009-07756 and
INCITE09207151PR (Spain).}



\maketitle

\section{Introduction}

A proper isometric action of a Lie group on a Riemannian manifold is called
{\em polar} if there exists a connected embedded submanifold~$\Sigma$ which
meets all orbits and which intersects the group orbits orthogonally in each of
its points. Such a submanifold~$\Sigma$ is called a {\em section} of the group
action. In the special case where the section is flat in the induced metric,
the action is called {\em hyperpolar}.

In this note, we prove some infinitesimal criteria for isometric actions to be
polar and apply them to obtain a classification of polar actions on symmetric
spaces with a fixed point. We neither assume that the symmetric space acted
upon is compact nor that it is irreducible. The result is the following.

\begin{theorem}\label{thMain}
Let $M$ be a symmetric space. Let $H$ be a connected closed subgroup of the
isometry group of~$M$ acting on $M$ with a fixed point. Then the $H$-action on
$M$ lifts to an isometric action with a fixed point on the universal
cover~$\tilde M$ which is polar if and only if the $H$-action on~$M$ is. In
this case, the $H$-action on $\tilde M$ is orbit equivalent to an action
defined as follows. Let $\tilde M = M_0 \times M_1 \times \ldots \times
M_{\ell}$ be a decomposition of the symmetric space~$\tilde M$ into one
Euclidean and $\ell$~irreducible factors. Let $H_0, \ldots, H_{\ell}$ be
compact Lie groups such that each factor $H_i$ acts trivially or polarly with a
fixed point on $M_i$ as described below.
\begin{enumerate}

\item[(i)] If $M_i$ is Euclidean, then the $H_i$-action on $M_i$ is orbit
    equivalent to a polar representation.

\item[(ii)] If $M_i$ is irreducible of rank greater than one, then the
    action of $H_i$ is orbit equivalent to the isotropy action of $M_i$.

\item[(iii)] If $M_i$ is isometric to $\SS^n$ or $\HH^n$, the action is
    given by a polar representation on $\R^n$.

\item[(iv)] If $M_i$ is isometric to $\C\P^n$ or $\C\HH^n$, then the action
    is given by a polar action on $\C\P^{n-1}$, see
    Section~\ref{seRankOne}.

\item[(v)] If $M_i$ is isometric to $\H\P^n$ or $\H\HH^n$, then the action
    is given by a polar action on $\H\P^{n-1}$, see
    Section~\ref{seRankOne}.

\item[(vi)] If $M_i$ is isometric to $\Ca\P^2$ or $\Ca\HH^2$, then the
    action is given by the isotropy group $\Spin(9)$ of~$M_i$ or one of its
    subgroups $\Spin(8)$, $\Spin(7)\cdot\SO(2)$, $\Spin(6)\cdot\Spin(3)$.

\end{enumerate}
\end{theorem}

Note that an action as described in Theorem~\ref{thMain} is hyperpolar if and
only if the actions of~$H_i$ on~$M_i$ are orbit equivalent to the isotropy
action of~$M_i$ for $i = 1,\ldots,n$. Hyperpolar actions without fixed point on
compact irreducible symmetric spaces have been classified in~\cite{K02}.

We remark that the result was known for irreducible spaces of higher
rank~\cite{Brk98} and for the compact symmetric spaces of rank one~\cite{PT99}.
It follows from Corollary~\ref{coPolarWithFP} that the actions described in
parts~(i) to~(vi) of the Theorem exhaust all orbit equivalence classes of polar
actions with a fixed point on strongly isotropy irreducible Riemannian
homogeneous spaces.

Using the results of Bergmann~\cite{Bgm01} on reducible polar representations,
all proper polar actions with a fixed point for each symmetric space (not only
up to orbit equivalence) can be described, although there is no convenient way
of writing down a complete list of all these actions due to the nature of the
result.

Our result is of interest for several reasons. One new aspect is that we use
here the duality between symmetric spaces of the compact and those of the
non-compact type. In particular, we obtain new results on polar actions on the
non-compact symmetric spaces of rank one. It is an intriguing question whether
duality can be used to study polar actions without fixed points.

It follows immediately from the main result of Dadok~\cite{Da85} that a polar
action with a fixed point on a simply connected symmetric space is orbit
equivalent to the product of polar actions on the individual factors of the
decomposition of the symmetric space into Euclidean and irreducible factors.
However, it is known that such a splitting into a product of actions does not
hold in general, i.e.\ for actions without fixed points, see the following
examples.

\begin{examples}\rm\hfill
\begin{enumerate}

\item Consider the action of $H = \SO(n)$ on $\R^{2n}=\R^n \oplus \R^n$
    where $H$ acts on each summand by the standard representation. Then the
    principal isotropy subgroup is $\SO(n-2)$ and $H$ acts with
    cohomogeneity one on the product $\SS^{n-1} \times \SS^{n-1}$ of the
    unit spheres in each summand~$\R^n$. This action is hyperpolar.

\item Let $G$ be a semisimple compact Lie group with a biinvariant metric
    and let $K_1,K_2 \subset G$ be two closed subgroups such that $(G,K_1)$
    and $(G,K_2)$ are \emph{symmetric pairs}, i.e.\ there is an involutive
    automorphism~$\sigma$ of~$G$ such that $(G^\sigma)_0 \subseteq K_i
    \subseteq G^\sigma$. The action of $K_1 \times K_2$ on~$G$ given by
    \begin{equation*}
    (k_1,k_2) \cdot g = k_1\, g\, k_2^{-1}
    \end{equation*}
    is hyperpolar~\cite{HPTT94} (a so-called \emph{Hermann action}). By
    \cite[Proposition~2.11]{HPTT94}, it follows that the action of $K_1
    \times K_2 \times G$ on~$G \times G$ given by
    \begin{equation*}
    (k_1,k_2,g) \cdot (g_1,g_2) = (k_1\, g_1\, g^{-1}, k_2\, g_2\, g^{-1})
    \end{equation*}
    is hyperpolar as well. Using \cite[Proposition~2.11]{HPTT94} once more,
    we see that the action of $\Delta G = \{(g,g)\mid g\in G\}$ on $G/K_1
    \times G/K_2$ is hyperpolar, but not a product. Example~(1) is the
    special case of this construction where $G=\SO(n)$, $K_1=K_2=\SO(n-1)$.

\item The action of $H = \Sp(2) \cong \Spin(5)$ on $\SS^7 \times \SS^4$,
    induced by the direct sum of the standard $\Sp(2)$-representation on
    $\H^2 = \R^8$ and the standard $\SO(5)$-representation on $\R^5$ is of
    cohomogeneity one, but does not split.

\item Similarly, the action of $\Spin(9)$ on $\SS^8 \times \SS^{15}$ by the
    standard plus the spin representation is also of cohomogeneity one and
    non-splitting.

\item The action of $\Spin(8)$ on $\SS^7 \times \SS^7 \times \SS^7$ given
    by the three inequivalent representations on $\R^8$ of $\Spin(8)$ is of
    cohomogeneity one, but does not split.

\end{enumerate}
\end{examples}
It would be desirable to have structural and classification results on
(hyper)polar actions on reducible symmetric spaces. Using the construction in
Example~(2), one may produce non-splitting actions with arbitrarily many
factors.

After the completion of the classification of polar and coisotropic actions on
the compact irreducible Hermitian symmetric spaces in his article~\cite{Bil04},
which showed in particular that polar actions are hyperpolar on the irreducible
Hermitian spaces of \emph{higher rank}, i.e.\ of rank greater than one,
Biliotti stated the following conjecture:

\begin{conjecture}\cite{Bil04}
A polar action on an irreducible compact symmetric space of rank greater than
one is hyperpolar.
\end{conjecture}

The conjecture has been proven to hold for the symmetric spaces with simple
compact isometry group~\cite{K07} and, more recently, for the exceptional
simple compact Lie groups~\cite{K09} with biinvariant metric. (It is still open
for the classical compact Lie groups.) However, in the noncompact setting it is
not true that polar actions on irreducible symmetric spaces of higher rank are
hyperpolar, since there are polar homogeneous foliations with non-flat sections
by recent results of~\cite{BDRT08}. But our result shows the conjecture of
Biliotti remains true also for spaces of noncompact type if we restrict
ourselves to polar actions of compact Lie groups, since these actions will have
a fixed point by Cartan's fixed point theorem. On the other hand, there are no
nontrival polar homogeneous foliations on irreducible compact symmetric spaces:

\begin{proposition}
Let $N$ be a symmetric space of compact type with a homogeneous polar
foliation. Let $\Sigma$ be a section and let $M$ be a principal orbit of the
homogeneous foliation. Then $\tilde N$ is isometric to the Riemannian product
$\tilde M \times \tilde \Sigma$, where the tilde denotes universal covers.
\end{proposition}

\begin{proof}
By \cite[Proposition 4.1]{K07}, the polar homogeneous foliation on $N$ lifts to
a homogeneous polar foliation on the universal cover $\tilde N$. It follows
from \cite[Lemma 1A.3]{PT99} that the corresponding isometric action on $\tilde
N$ does not have any exceptional orbits. Since a homogeneous foliation does not
have singular orbits, it follows that the Weyl group $W_{\Sigma}$ is trivial,
cf.\ \cite[Lemma 5.1]{K07}. Therefore, one may apply the Splitting Theorem~5.2
in \cite{K07} and it follows that $\tilde N$ is isometric to  $\tilde M \times
\tilde \Sigma$.
\end{proof}

It is an interesting open question whether there are polar actions with
non-flat section and singular orbits on irreducible symmetric spaces of
non-compact type.

We briefly summarize the contents of this paper. In Section~\ref{seCriteria} we
prove a criterion for polarity of actions with a fixed point. In
Section~\ref{seDuality} we use this criterion to show that the duality of
symmetric spaces leads to a one-to-one correspondence between polar actions
with a fixed point on a pair of dual symmetric spaces. We classify polar
actions with a fixed point in rank one symmetric spaces in
Section~\ref{seRankOne}. Theorem~\ref{thMain} is proven in
Section~\ref{seProof}.

\section{Criteria for polarity}\label{seCriteria}

We say that two isometric Lie group actions on two Riemannian manifolds $M_1$
and $M_2$, respectively, are \emph{orbit equivalent}, if there is an isometry
$M_1 \to M_2$ which maps connected components of orbits onto connected
components of orbits. Two isometric actions of a Lie group~$H$ on two
Riemannian manifolds $M_1$ and $M_2$ are called \emph{conjugate} if there is an
isometry $F \colon M_1 \to M_2$ such that $F(h \cdot p) = h \cdot F(p)$ for all
$h \in H$ and $p \in M_1$. Obviously, the polarity of an action only depends on
the orbit equivalence class of the action.

We will now prove an infinitesimal criterion which allows to decide if the
orbits of an isometric Lie group action intersect a totally geodesic
submanifold orthogonally. This is a simple but fundamental observation from
which one may deduce for instance the criteria for polarity of actions on
symmetric spaces which had been obtained in various special cases by different
authors, cf.~\cite[Theorem~2.1 and Corollary~2.12]{HPTT94},
\cite[Proposition]{G04}, \cite[Proposition~4.1]{K07},
\cite[Theorem~4.1]{BDRT08}.

\begin{lemma}\label{leOrthKill}
Let $M$ be a Riemannian manifold and let $\Sigma$ be a connected totally
geodesic submanifold of~$M$. Let $p\in\Sigma$ and let $X$ be a Killing vector
field. Then $X(q) \in N_q\Sigma$ holds for all $q \in \Sigma$ if and only if
$X(p) \in N_p\Sigma$ and $\nabla_v X \in N_p \Sigma$ for all $v \in T_p\Sigma$.
\end{lemma}

\begin{proof}
Let $\gamma \colon \R \to \Sigma$ be a geodesic such that $\gamma(0)=p$. Then
$J(t) := X(\gamma(t))$ is a Jacobi field along~$\gamma$. Let $n = \dim M$. Let
$e_1(t), \ldots, e_n(t)$ be parallel orthonormal fields along $\gamma$ such
that $e_1(0), \ldots, e_k(0) \in T_p\Sigma$ for $k = \dim \Sigma$. Since
$\Sigma$ is totally geodesic, we have that the vectors $e_1(t),\ldots,e_k(t)$
are contained in $T_{\gamma(t)}\Sigma$ for all $t$.  Let $J(t) = \sum_i f_i(t)
e_i(t)$ for some functions $f_i$. The Jacobi equation $\frac{D^2J}{dt^2} = R(
\gamma^\prime(t),J(t))\gamma^\prime(t)$ is then equivalent to the system
\begin{equation*}
f^{\prime\prime}_j(t) + \sum_i a_{ij}(t)f_i(t) = 0, \quad j=1,\ldots,n,
\end{equation*}
where the functions $a_{ij}$ are given by $a_{ij}(t) = \langle
R(\gamma^\prime(t),e_i(t))\gamma^\prime(t),e_j(t)\rangle$, see~\cite{dC92}. It
follows from the symmetries of the Riemann tensor that the matrix $a_{ij}(t)$
is symmetric for all~$t$. Moreover, since $\gamma^\prime(t), e_1(t),\ldots,
e_k(t)$ are tangent to the totally geodesic submanifold $\Sigma$, we have
$a_{ij}(t) = a_{ji}(t) = 0$ if $i \le k$ and $j > k$. It follows that $f_1(t) =
f_1^\prime(t) = \ldots = f_k(t) = f_k^\prime(t) = 0$ for all~$t$ (which is
equivalent to $J(t) \in N_{\gamma(t)}\Sigma$ for all~$t$) if and only if we
have $f_1(0) = f_1^\prime(0) = \ldots = f_k(0) = f_k^\prime(0) = 0$ (which is
equivalent to $J(0)\in N_p\Sigma$ and $J^\prime(0) = \frac{\nabla J}{dt}(0) =
\nabla_{\gamma^\prime(0)} X \in N_p\Sigma$). This shows the statement of the
lemma.
\end{proof}

From the above lemma, we obtain a general criterion for polarity of isometric
actions.

\begin{corollary}\label{coPolarCrit}
Let $M$ be a complete connected Riemannian manifold and let $\Sigma$ be a
connected totally geodesic embedded submanifold of~$M$. Then the proper
isometric action of a Lie group $H$ on~$M$ is polar with section~$\Sigma$ if
and only if there is a point $p \in \Sigma$ such that (i)~$T_p\Sigma \subseteq
N_p(H \cdot p)$, (ii)~the slice representation of the isotropy subgroup $H_p$
on $N_p(H \cdot p)$ is such that $T_p\Sigma$ meets all $H_p$-orbits and
(iii)~$\nabla_vX \in N_p\Sigma$ for all $v \in T_p\Sigma$ and all Killing
vector fields $X$ induced by the $H$-action on~$M$.
\end{corollary}

\begin{proof}
Let $M$ be a Riemannian $H$-manifold and let $\Sigma$ be a connected totally
geodesic embedded submanifold such that (i) and (ii) hold for some point $p \in
\Sigma$. We show that then $\Sigma$ intersects all $H$-orbits. Let $q \in M$.
Since the orbit $H \cdot p$ is closed~\cite{DR08} and $M$ is complete and
connected, there is a shortest geodesic~$\gamma$ joining $q$ and $H \cdot p$,
which intersects $H \cdot p$ orthogonally at some point $h \cdot p$. Then
$h^{-1} \circ \gamma$ is a geodesic joining $p$ and $h^{-1}(q)$. Since the
slice representation of~$H_p$ on the normal space $N_p (H \cdot p)$ is such
that $T_p\Sigma$ meets all orbits, there is an element $g \in H_p$ such that
$g^{-1} \circ h^{-1} \circ \gamma$ is a geodesic joining the two points $p$ and
$g^{-1}(h^{-1}(q))$ whose image is contained in~$\Sigma$. This shows that the
$H$-orbit through~$q$ intersects~$\Sigma$. The rest of the statement follows
immediately from Lemma~\ref{leOrthKill}.
\end{proof}

The key observation for this article is the following criterion for polarity of
isometric actions with a fixed point, which enables us to use the duality of
symmetric spaces.

\begin{theorem}\label{thCriterium}
Let $M$ be a connected Riemannian manifold with isometry group~$I(M)$ and let
$p \in M$. A closed subgroup $H \subseteq I(M)_p$ acts polarly on~$M$ if and
only if the slice representation of~$H = H_p$ on~$T_pM$ is polar and for any of
its sections~$\g{s}$, the exponential image $\Sigma := \exp_p(\g{s})$ is an
embedded totally geodesic submanifold of~$M$. In this case, $\Sigma$ is a
section for the $H$-action on~$M$.
\end{theorem}

\begin{proof}
Assume $H$ acts polarly on~$T_pM$ and for any of its sections~$\g{s}$, the
exponential image $\exp_p(\g{s})$ is an embedded totally geodesic submanifold
of~$M$. Let $x$ be an element of the Lie algebra of~$H$. Then for all $q \in
M$, the Killing vector field~$X$ corresponding to~$x$ is given by $X(p) =
\left. \frac d {ds} \right|_{s=0} \left( h_s(p) \right),$ where $h_s$ denotes
the isometry of $M$ given by the group element $\exp(sx)$. Let $\exp_p \colon
T_pM \to M$ denote the Riemannian exponential map of~$M$ at the point~$p$ and
let $v \in \g{s}$. Then we have
\begin{align*}
\nabla_vX &=
 \left. \frac\nabla{\partial t} \frac\partial{\partial s} h_s(\exp_p(tv))
\right\rvert_{s=t=0} = \left. \frac\nabla{\partial s} \frac\partial{\partial t}
h_s(\exp_p(tv)) \right\rvert_{s=t=0} =\\ &= \left. \frac\nabla{\partial s}\left(
\frac\partial{\partial t} \left.h_s(\exp_p(tv))\right\rvert_{t=0} \right)
\right\rvert_{s=0} = \frac\nabla{\partial s} \left. (h_s)_{*p}(v)
\right\rvert_{s=0} \in N_p\Sigma,
\end{align*}
since the slice representation of~$H$ on~$T_pM$ is polar with section~$\g{s}$.
Now it follows from Corollary~\ref{coPolarCrit} that the $H$-action is polar
with section~$\Sigma$. The converse is a well known fact.
\end{proof}

To find all polar actions on a given space, it is often useful to know \emph{a
priori} that certain subgroups of a group acting isometrically do not act
polarly. The following theorem gives such an information for irreducible
representations which are not transitive on the sphere.

\begin{theorem}\label{thPolarSubaction}\cite{KP03}
Let $G\subset \SO(n)$ be a closed connected subgroup which acts irreducibly on
$\R^n$ and non-transitively on the sphere $\SS^{n-1}\subset\R^n$. Let $H\subset
G$ be a closed connected subgroup $\neq\{e\}$ that acts polarly on $\R^n$. Then
the $H$-action and the $G$-action on $\R^n$ are orbit equivalent.
\end{theorem}

This can be used to show the following fact on polar actions with a zero- or
one-dimensional orbit on strongly isotropy irreducible Riemannian homogeneous
spaces. In particular, it follows that polar actions with a fixed point on
irreducible symmetric spaces of rank greater than one are hyperpolar and in
fact orbit equivalent to the isotropy action. This had been shown already
in~\cite[Theorem~2.2]{Brk98}.

\begin{corollary}\cite[Corollary~6.2]{K07}\label{coPolarWithFP}
Let $X$ be a strongly isotropy irreducible Riemannian homogeneous space. Assume
a connected compact Lie group $H$ acts polarly and non-trivially on~$X$ and the
$H$-action has a one-dimensional orbit $H \cdot p$ or a fixed point $p \in X$.
Then the space~$X$ is locally symmetric. Furthermore, $X$ is a rank-one
symmetric space or the action of $H$ is orbit equivalent to the action of the
connected component of the isotropy group of~$X$ at~$p$.
\end{corollary}

\section{Duality}\label{seDuality}

The statement of Theorem~\ref{thMain} was previously known in the case of
irreducible symmetric spaces of higher rank~\cite{Brk98} and for compact rank
one symmetric spaces it follows from~\cite{PT99}.  Thus it essentially remains
to extend the classification to the noncompact rank-one symmetric spaces.
Theorem~\ref{thCriterium} enables us to do this by using duality between
symmetric spaces of the compact and the non-compact type.

Let $(G,K)$ be a symmetric pair and let $\g{g} = \g{k} \oplus \g{p}$ be an
$\Ad(K)$-invariant decomposition. Then one may define a Lie algebra $\g{g}^*$
be setting $\g{g}^* = \g{k} \oplus i\g{p}$, where $i = \sqrt{-1}$. If a
symmetric pair $(G^*,K)$ is such that $\g{g}^*$ is isomorphic to the Lie
algebra of~$G^*$, we say that $(G^*,K)$ is \emph{dual} to $(G,K)$.

\begin{theorem}\label{thDualPolar}
Let $(G,K)$ be a symmetric pair such that $G/K$ is an irreducible space of the
compact or noncompact type and let $(G^*,K)$ be a dual symmetric pair. Let $H
\subseteq K$ be a closed subgroup. Then the action of $H$ on $G/K$ is polar if
and only if the action of $H$ on $G^*/K$ is polar.
\end{theorem}

\begin{proof}
Assume first that $(G,K)$ is a symmetric pair and $G/K$ is of rank one. Let
$(G^*,K)$ be a dual of $(G,K)$. Let $H$ be a closed subgroup of~$K$ acting
polarly on $G/K$ with a fixed point. By Theorem \ref{thCriterium}, $H$ acts
polarly on $T_{eK}(G/K)$. Let $\Sigma\subset T_{eK}(G/K)$ be a section of this
action. Then $\Sigma$ is a Lie triple system in $\g{p}=T_{eK}(G/K)$ and hence
$i\Sigma$ is a Lie triple system in $i\g{p}=T_{eK}(G^*/K)$. In a rank one
symmetric space, every totally geodesic submanifold is embedded, hence it
follows from Theorem \ref{thCriterium} that the action of~$H$ on $G^*/K$ is
polar if $\rank G^*/K =\rank G/K =1$. For spaces of higher rank, the assertion
follows immediately from Corollary~\ref{coPolarWithFP}.
\end{proof}

\begin{corollary}\label{coSectConstCurv}
Let $\Sigma$ be the section of a polar action with fixed point on an
irreducible symmetric space of noncompact type. Then $\Sigma$ is isometric to a
product $\R^{n_0} \times \HH^{n_1} \times \ldots \times \HH^{n_k}$.
\end{corollary}

\begin{proof}
This follows via Theorem~\ref{thDualPolar} from an analogous statement for
spaces of compact type \cite[Theorem~5.4]{K07}.
\end{proof}

We remark that Corollary~\ref{coSectConstCurv} does not hold anymore if one
drops the assumption that the polar action has a fixed point, see
\cite[Proposition 4.2]{BDRT08}.

\section{Actions on rank one symmetric spaces}\label{seRankOne}

Let us recall some results of~\cite{PT99} about polar actions on compact rank
one symmetric spaces. We say that an isometric action of a connected compact
Lie group $H \subseteq \U(n+1)$ on $\C\P^n = \U(n+1) / (\U(n) \times \U(1))$ is
\emph{induced by the isotropy representation of a Hermitian symmetric space}
$U/K$ if there is a semisimple Hermitian symmetric space~$U/K$ such that
$(U,K)$ is a symmetric pair and the natural action of $H$ on $\C^{n+1}$ is
equivalent to the isotropy representation of $U/K$.

Similarly, let $H_i \subseteq \Sp(n_i)$ be connected closed subgroups for $i =
1,\ldots,k$; we say that the isometric action of $H_1 \times \ldots \times H_k
\subseteq \Sp(n_1)\times \ldots \times \Sp(n_k)\subseteq \Sp(n+1)$ on $\H\P^n$,
where $n+1 = n_1 +\ldots+ n_k$, is \emph{induced by the isotropy representation
of a product of $k$~quaternion-K\"{a}hler symmetric spaces} if there are
quaternion-K\"{a}hler symmetric spaces $U_1/K_1, \ldots, U_k/K_k$ such that the
action of~$H_i \times \Sp(1)$ on~$\H^{n_i}$ (given by restriction of the usual
$\Sp(n_i) \times \Sp(1)$-representation) is equivalent (on the Lie algebra
level) to the isotropy representation of $U_i/K_i$.

It was proved in~\cite{PT99} that the actions on $\C\P^n$ as described above
are polar and exhaust the orbit equivalence classes of polar actions
on~$\C\P^n$. Furthermore, the actions on $\H\P^n$ where at most one of the
quaternion-K\"{a}hler symmetric spaces $U_i/K_i$ is of rank~$\ge 2$, are polar and
exhaust the orbit equivalence classes of polar actions on~$\H\P^n$.

We remark that an action on $\C\P^n$, resp.\ $\H\P^n$, induced by the isotropy
representation of a Hermitian symmetric space, resp.\ a product of
$k$~quaternion-K\"{a}hler symmetric spaces, has fixed point if and only if the
symmetric space decomposes as a Riemannian product containing a $\C\P^1$,
resp.\ an $\H\P^1$, factor.

The next proposition shows that the orbit equivalence classes of polar actions
with a fixed point on ${\mathbb K}\P^n$ resp.\ ${\mathbb K}\HH^n$ are in
one-to-one correspondence with the orbit equivalence classes of polar actions
on ${\mathbb K}\P^{n-1}$ for ${\mathbb K}=\R, \C, \H$.

\begin{proposition}\label{prRank1Polar}
Let $n \ge 2$.

\begin{enumerate}
\item[(i)] Let $H$ be closed subgroup of~$\SO(n)$. Then the group $H$ acts
    polarly on $\SS^n = \SO(n+1)/\SO(n)$ or $\HH ^n = \SO(n,1)_0/\SO(n)$,
    respectively, if and only if $H$ acts polarly on~$\SS^{n-1}$ by
    restriction of the $n$-dimensional standard representation of $\SO(n)$.

\item[(ii)] Let $H$ be a closed subgroup of~$\U(n)$. If the action of~$H$
    on $\C\P^n=\SU(n+1)/\U(n)$ or $\C\HH^n=\SU(n,1)/\U(n)$, respectively,
    is polar, then the $H$-action on $\C\P^{n-1} = \U(n)/(\U(n-1)\times
    \U(1))$ is polar. Conversely, if the action of $H$ on $\C\P^{n-1}$ is
    polar and nontrivial, then the action of $\U(1) \cdot H \subseteq
    \U(n)$ on $\C\P^n$ or $\C\HH^n$, respectively, is polar with a fixed
    point.

\item[(iii)] Let $H$ be a closed subgroup of~$\Sp(n) \times \Sp(1)$. If the
    action of~$H$ on $\H\P^n=\Sp(n+1)/(\Sp(n)\times\Sp(1))$ or
    $\H\HH^n=\Sp(n,1)/(\Sp(n)\times\Sp(1))$, respectively, is polar, then
    the $\pi(H)$-action on $\H\P^{n-1} = \Sp(n)/(\Sp(n-1) \times \Sp(1))$,
    where $\pi \colon \Sp(n) \times \Sp(1) \to \Sp(n)$ is the natural
    projection onto the first factor, is polar. Conversely, let $H$ be a
    closed subgroup of $\Sp(n)$ acting polarly and nontrivially on
    $\H\P^{n-1}$, then $H \times \Sp(1)$ acts polarly with a fixed point on
    $\H\P^n$ or $\H\HH^n$, respectively.

\item[(iv)] Let $H$ be a closed connected subgroup of~$\Spin(9)$. Then the
    action of $H$ on $\Ca\P^2 = \LF_{4}  / \Spin(9)$ or $\Ca\HH^2 =
    \LF_{4}^* / \Spin(9)$, respectively, is polar if and only if $H$ is
    conjugate to one of the following: $\Spin(9)$, $\Spin(8)$,
    $\Spin(7)\cdot\SO(2)$, $\Spin(6)\cdot\Spin(3)$.

\end{enumerate}
\end{proposition}
\begin{proof}
In view of Theorem~\ref{thDualPolar}, it suffices to prove the proposition in
the compact case. Part (i) follows immediately from the fact that polar actions
on the sphere $\SS^n$ are precisely the restrictions of polar representations
on~$\R^{n+1}$.

We will now prove part (ii).  Let $H \subseteq \U(n)$ be a closed subgroup
acting polarly (and with a fixed point) on $\C\P^n$, homogeneously presented as
$\SU(n+1)/\U(n)$, where $\U(n)$ is embedded into $\SU(n+1)$ as the subgroup
$\SUxU n1$. The action of this group leaves invariant a totally geodesic
$\C\P^{n-1}$, which we may identify with $\U(n) / (\U(n-1) \times \U(1))$. By
\cite[Lemma 4.2]{K07}, the restriction of the $H$-action to this submanifold is
polar.

Now assume $H \subseteq \U(n)$ is a closed subgroup acting polarly and
nontrivially on $\C\P^{n-1} = \U(n) / (\U(n-1) \times \U(1))$. By the results
of \cite{PT99}, the $H$-action on $\C\P^{n-1}$ is orbit equivalent to an action
induced by the isotropy representation of a Hermitian symmetric space and we
may assume that the action of~$H$ on~$\C^n$ is equivalent to the isotropy
representation of a Hermitian symmetric space~$X$. Then the action of $\U(1)
\cdot H \subseteq \U(n)$ on $\C\P^n = \SU(n+1) / \U(n)$ is orbit equivalent to
the action induced by the Hermitian symmetric space $X \times \C\P^1$, showing
that the action of $\U(1) \cdot H$ on $\C\P^n$ is polar \cite{PT99} with a
fixed point.

To prove part (iii), assume $H \subset \Sp(n) \times \Sp(1)$ is a closed
subgroup acting polarly (and with a fixed point) on $\H\P^n$. The action of
this group leaves invariant a totally geodesic $\H\P^{n-1}$. By \cite[Lemma
4.2]{K07}, the restriction of the $H$-action to this submanifold is polar.

Now assume $H \subseteq \Sp(n)$ is a closed subgroup acting polarly and
nontrivially on $\H\P^{n-1}$. By the results of \cite{PT99}, the $H$-action on
$\H\P^{n-1}$ is orbit equivalent to an action induced by the isotropy
representation of a product of $k$~quaternion-K\"{a}hler symmetric spaces $U_1/K_1
\times \ldots \times U_k/K_k$ where at most one factor is of higher rank. Then
the action induced by the isotropy representation of $U_1/K_1 \times \ldots
\times U_k/K_k \times \H\P^1$ on $\H\P^n$ is orbit equivalent to the action of
$H \times \Sp(1)$, which is hence polar (and has a fixed point).

Part~(iv) follows immediately from~\cite{PT99}.
\end{proof}

\section{Proof of Theorem~\ref{thMain}}\label{seProof}

\begin{proposition}\label{prSplitting}
Let $(G_0,K_0), \ldots, (G_{\ell},K_{\ell})$ be symmetric pairs such that
$G_0/K_0$ is of Euclidean type and $G_1/K_1, \ldots, G_{\ell}/K_{\ell}$ are
irreducible symmetric spaces. Define $G = G_0 \times G_1 \times \ldots \times
G_{\ell}$ and $K = K_0 \times K_1 \times \ldots \times K_{\ell}$. Let $H
\subseteq K$ be a closed subgroup. Then there is a natural action of~$H$ on the
spaces $G_i / K_i$ given by $\pi_i \circ \iota$, where $\iota \colon H \to K$
is the inclusion and $\pi_i \colon K \to K_i$ is the projection. If the action
of~$H$ on $G/K$ is polar then the action of~$H$ on each of the spaces $G_i /
K_i$ is polar. Moreover, the $H$-action on $G/K$ is orbit equivalent to the
product action of $\Pi_{i=1}^\ell H$ on $G_1/K_1 \times \ldots \times
G_{\ell}/K_{\ell}$.
\end{proposition}

\begin{proof}
Assume the $H$-action on $G/K$ is polar. Then by Theorem~\ref{thCriterium}, the
slice representation of~$H$ on $T_{eK}G/K$ is polar as well and for any of its
sections $\g{s}$, the exponential image $\exp_{eK}(\g{s})$ is an embedded
totally geodesic submanifold of~$G/K$. Since $H \subseteq K$, the group $H$
leaves each $T_{eK}G_i/K_i$ invariant. Then \cite[Theorem~4]{Da85} implies that
the $H$-action on each invariant summand $T_{eK}G_i/K_i$ is also polar and the
section $\g{s}$ is of the form $\g{s}_0 \oplus \ldots \oplus \g{s}_{\ell}$,
where $\g{s}_i$ is a section of the $H$-action on $T_{eK}G_i/K_i$. From this we
see that $\exp_{eK_i}(\g{s}_i)$ is an embedded totally geodesic submanifold of
$G_i/K_i$ for each $i=0,\ldots,\ell$. By using Theorem~\ref{thCriterium} once
more, we get that the $H$-action on each factor $G_i/K_i$ is polar with
section~$\exp_{eK_i}(\g{s}_i)$. The last statement follows since the orbits of
the $H$-action on $G/K$ are certainly contained in those of the $\Pi_{i=1}^\ell
H$-action and both actions are of the same cohomogeneity.
\end{proof}

\begin{proof}[Proof of Theorem~\ref{thMain}]
Let $M = G/K$ be an effective presentation and $G = \Isom(M)$. We may assume
that $H$ is a closed subgroup of~$K$. Let $\tilde G$ be the universal cover
of~$G$ with covering map $\pi \colon \tilde G \to G$. Define $\tilde K :=
\pi^{-1}(K)$ and identify $M = \tilde G / \tilde K$. Then $\pi|_{\tilde K}
\colon \tilde K \to K$ is a covering map whose kernel agrees with the
effectivity kernel of the $\tilde K$-action on $M = \tilde G / \tilde K$. Let
$h \in H$. Since $H$ is connected, there is a path~$h(t)$ defined on~$[0,1]$
and connecting $h = h(1)$ with the identity element $e = h(0)$. Let $\tilde
h(t)$ be the continuous lift of $h(t)$ such that $\pi(\tilde h(t))=h(t)$ and
$\tilde h(0) = e$. Then $\tilde h := \tilde h(1)$ lies in the connected
component $\tilde K_0$. Define an action of~$H$ on $\tilde M$ by requiring that
an element $h \in H$ acts on $\tilde M = \tilde G / \tilde K_0$ by $\tilde h$
as constructed above. Since the action of the isometry $\tilde h$ is determined
by its action on $T_{e\tilde K} \tilde M$ and the representation of $\tilde
K_0$ on $T_{e\tilde K} \tilde M$ is equivalent to $\rho \circ \pi \colon \tilde
K_0 \to {\rm O}(T_{eK} M)$, where $\rho$ is the isotropy representation of~$K$
on~$M$, the action of $H$ on $\tilde M$ is well-defined, i.e.\ does not depend
on the path $h(t)$.  It is now clear that we have defined an isometric Lie
group action of~$H$ on~$\tilde M$.

It follows from Theorem~6 that the action of~$H$ on $\tilde M$ is polar if and
only if the action of~$H$ on~$M$ is. Assume the $H$-action on $\tilde M$ is
polar. It follows from Proposition~\ref{prSplitting} that the induced actions
on the factors $M_i$ are all polar, hence either trivial, as described in
Proposition~\ref{prRank1Polar} or Corollary~\ref{coPolarWithFP}, or given by a
polar linear representation if $M_i$ is Euclidean.
\end{proof}

\end{document}